\documentclass{article}
\usepackage{amssymb}
\begin{document}
\input psbox.tex
\psfordvips

\centerline{\LARGE\bf On the classification of convex lattice}

\medskip
\centerline{\LARGE\bf polytopes }

\bigskip \centerline{\bf Heling Liu and Chuanming Zong\footnote{The second author
is supported by 973 Program 2011CB302400,
the National Natural Science Foundation of China (Grant No.
11071003), the Chang Jiang Scholars Program of China and LMAM at
Peking University.}}

\bigskip
\centerline{\it Dedicated to Professor Peter M. Gruber on the
occasion of his 70th birthday}

\bigskip
\noindent {\bf Abstract.} In 1980, Arnold studied the classification problem
for convex lattice polygons of given area. Since then this problem and its
high dimensional analogue have been studied by B\'ar\'any, Pach, Vershik and others.
Bounds for the number of non-equivalent $d$-dimensional convex lattice polytopes
of given volume have been achieved. In this paper we study Arnold's problem for
centrally symmetric lattice polygons and the classification problem for convex
lattice polytopes of given cardinality. In the plane we obtain analogues to the bounds
of Arnold, B\'ar\'any and Pach in both cases. However, the number of non-equivalent
$d$-dimensional convex lattice polytopes of $w$ lattice points is infinite whenever
$w-1\ge d\ge 3$, which may intuitively contradict to B\'ar\'any and Vershik's
upper bound.


\section*{1. Introduction}

A convex lattice polytope in $\mathbb{E}^d$ is the convex hull of a finite subset of
the integer lattice $\mathbb{Z}^d$. Equivalently, it is a convex polytope, all vertices
of which are in $\mathbb{Z}^d$. One of the first results about lattice
polygons was discovered by G. Pick in 1899, which is known as Pick's
theorem. In 1967, E. Ehrhart discovered his polynomiality theorem.
In recent years, Ehrhart's polynomials have attracted the attention of
many authors (see \cite{barv97}, \cite{grit93} or \cite{grub07}). In
fact, up to now the knowledge about convex lattice polytopes is
still very limited (see \cite{gard01} and \cite{gard05}).

Let $P$ denote a $d$-dimensional convex lattice polytope, let $v(P)$ denote the volume
of $P$, and let $|P|$ denote the cardinality of $P\cap \mathbb{Z}^d$.

Let $P_1$ and $P_2$ be two $d$-dimensional lattice polytopes. If
there is a $\mathbb{Z}^d$-preserving affine transformation $\sigma$ satisfying
$$P_2=\sigma (P_1),$$
then we say $P_1$ and $P_2$ are equivalent. For convenience,
we write $P_1 \sim P_2$ for short. It is easy to see that, if
$P_1\sim P_2$ and $P_2\sim P_3$, then we have $P_1\sim P_3$. In
addition, if $P_1\sim P_2$, then we have
$$v(P_1)=v(P_2)$$
and
$$|P_1|=|P_2|.$$

Clearly, the equivalence relation $\sim$ divides convex lattice
polytopes into different classes. Using triangulations, it can be easily shown that
$$d!\cdot v(P)\in \mathbb{Z}$$
holds for any $d$-dimensional convex lattice polytope $P$. Let $v(d,m)$ denote the number of
different classes of the $d$-dimensional convex lattice polytopes $P$
with $v(P)=m/d!$, where both $d$ and $m$ are positive integers. In
1980, Arnold \cite{arno80} studied the values of $v(2,m)$ and proved
$$m^{1\over 3}\ll \log v(2,m)\ll m^{1\over 3}\log m.\eqno(1)$$

\medskip\noindent
{\bf Remark 1.} In this paper $f(d,m)\ll g(d, m)$ means that, for fixed positive integer
$d$,
$$f(d, m)\le c_d\cdot g(d,m)$$
holds for all positive integers $m$ with a suitable constant $c_d$.

\bigskip
In 1992, B\'ara\'ny and Pach \cite{bara92} improved Arnold's
upper bound by removing the $\log m$ term; B\'ar\'any and Vershik \cite{bara92'}
obtained a general upper bound
$$ \log v(d,m)\ll m^{{d-1}\over {d+1}}.\eqno(2)$$
In the literature, some citations on $v(d,m)$ are confusing. We will clarify the situation in
Section 3.

Let $v^\ast(d,m)$ denote the number of different classes of the $d$-dimensional
centrally symmetric convex lattice polytopes $P$ with $v(P)=m/d!$, let $\kappa(d,w)$
denote the number of different classes of $d$-dimensional convex lattice polytopes
$P$ with $|P|=w$, and let $\kappa^\ast(d,w)$ denote the number of different classes of
$d$-dimensional centrally symmetric convex lattice polytopes $P$ with $|P|=w$. Then
we have $v^\ast(d,m)=0$ whenever $m$ is odd and $\kappa^\ast(d,w)=0$ if $w$ is even. Therefore
in this paper we assume that the $m$ in $v^\ast(d,m)$ is even and the $w$ in
$\kappa^\ast(d,w)$ is odd.

In this paper, we study Arnold's problem for the centrally symmetric lattice polygons
and the classification problem for convex lattice polytopes of given cardinality.
In Section 2 we introduce a basic lemma on the structures of convex lattice polytopes. In
Section 3 we review the known results about $v(d,m)$ and prove
$$m^{1\over 3}\ll \log v^\ast(2, m)\ll m^{1\over 3}.$$
In Section 4 we prove
$$w^{1\over 3}\ll \log \kappa(2, w)\ll w^{1\over 3}$$
and
$$w^{1\over 3}\ll \log \kappa^\ast(2, w)\ll w^{1\over 3}.$$
In Section 5 we show that
$$\kappa(d,w)=\infty $$ whenever $d\ge 3$ and $w\ge d+1$, which may intuitively contradict to (2)
about $v(d,m)$, and
$$\log \kappa^\ast (d, w)\ll w^{{d-1}\over {d+1}}.$$

\section*{2. Rabinowitz's Lemma}

In this section we introduce a basic result about the structures of
convex lattice polytopes which will be useful in Section 4. The result
was discovered by S. Rabinowitz in 1989 and was published at Utilitas Mathematica.
Since the result is elegant and the journal is hard to find, we reproduce its proof
here.

\medskip
\noindent {\bf Lemma 1 (\cite{rabi89}).} {\it Let $P$ be a $d$-dimensional convex
lattice polytope and let $m$ be a natural number satisfying $|P|\ge m^d+1$.
Then $P$ has at least $m+1$ collinear lattice points.}

\medskip\noindent
{\bf Proof.} Consider the coordinates of the integer points modulo $m$.
Since there are only $m^d$ distinct $d$-tuples of integers modulo
$m$, some two points ${\bf x}$ and ${\bf y}$ of $P\cap \mathbb{Z}^d$ must be
congruent (mod $m$). In other words, for all $i=1, 2, \ldots , d$, we have
$$x_i-y_i\equiv 0\quad ({\rm mod}\ m).$$
By convexity, all the $m+1$ collinear lattice points
$${\bf x}+\mbox{${j\over m}$}({\bf y}-{\bf x}),\quad j=0, 1, 2, \ldots
, m$$ belong to $P$. The lemma is proved. \hfill{$\square $}

\section*{3. Arnold's Problem}

In this section, we review the known results about $v(d,m)$ and prove
$$m^{1\over 3}\ll \log v^*(2,m)\ll m^{1\over 3}.$$

Let $\varepsilon (P)$ denote the cardinality of the vertices of $P$. To
prove the upper bound in (1), Arnold \cite{arno80} showed that
$$\varepsilon (P)\ll v(P)^{1\over 3}\eqno(3)$$
holds for all two-dimensional convex lattice polygons. In 1984, Konyagin and
Sevastyanov \cite{kony84} generalized (3) to $d$
dimensions by proving
$$\varepsilon (P)\ll v(P)^{{d-1}\over {d+1}}.\eqno(4)$$
In fact, this upper bound was first achieved by Andrews
\cite{andr65} in 1965.

At the end of \cite{arno80}, Arnold made a remark that \lq\lq In
$\mathbb{Z}^d$, $1/3$ is probably replaced by $(d-1)/(d+1)$. Proof of the
lower bound: let $x_1^2+\ldots +x_{d-1}^2\le x_d\le A$." Therefore,
the following problem is cited as Arnold's question in the
literature (see \cite{bara08} and \cite{bara08'}): {\it To
investigate $v(d,m)$ and to determine the order of magnitude of
$\log v(d,m)$.}

In 1992, B\'ar\'any and Pach \cite{bara92} improved Arnold's
upper bound to
$$\log v(2,m)\ll m^{1\over 3};\eqno(5)$$
B\'ar\'any and Vershik \cite{bara92'} generalized (5) to
$d$ dimensions by proving
$$\log v(d,m)\ll m^{{d-1}\over {d+1}}.\eqno(6)$$
In \cite{bara92'}, the authors attributed
$$\log v(d,m)\gg m^{{d-1}\over {d+1}}\eqno(7)$$
and $$\log v(d,m)\ll m^{{d-1}\over {d+1}}\log m$$ to \cite{arno80}
and \cite{kony84}, respectively. In fact, neither of them contains
such proofs. In particular, a proof for (7) seems non-trivial.
Therefore, to determine the order of magnitude of $\log v(d,m)$ for
fixed $d$ and large $m$ is still a basic open problem.

In \cite{bara08} and \cite{bara08'}, it was concluded that
$$\log \left(\sum_{j=1}^mv(d,j)\right)\gg
m^{{d-1}\over {d+1}}$$ and attributed this lower bound to Arnold
\cite{arno80}. Unfortunately, a rigorous proof is missing as well.

To estimate $v^\ast (2, m)$ we have the following result.

\medskip\noindent
{\bf Theorem 1.} {\it When $m$ is even and sufficiently large, we have}
$$m^{1\over 3}\ll \log v^\ast(2, m)\ll m^{1\over 3}.$$

To prove this theorem, we need the following technical lemma.

\medskip\noindent
{\bf Lemma 2.} {\it Let $T_\ell$ denote the lattice triangle with vertices ${\bf o},$
$(\ell,0)$ and $(0,\ell)$, let $S_\ell$ denote the lattice square with vertices ${\bf o},$
$(\ell,0)$, $(\ell,\ell)$ and $(0,\ell)$, and let $k$ be an integer with $\ell\le k\le \ell^2$.
Then, there is a convex lattice polygon $P$ satisfying both $T_\ell\subset P\subseteq S_\ell$
and} $$v(P)=\mbox{${1\over 2}$}(\ell^2+k).$$

\noindent
{\bf Proof.} Let $P_j$ denote the convex lattice pentagon with vertices ${\bf o}$, $(0,\ell)$,
$(1,\ell)$, $(1+j,$ $\ell-j)$ and $(\ell,0)$, let $H_{i,j}$ denote the hexagon with vertices ${\bf o}$,
$(0,\ell)$, $(i+1,\ell)$, $(2+i+j,\ell-j)$, $(\ell, i)$ and $(\ell,0)$, and let $H'_{i,j}$ denote the
hexagon with vertices ${\bf o}$, $(0,\ell)$, $(i+1,\ell)$, $(\ell-j,2+i+j)$, $(\ell,2+i)$ and $(\ell,0)$. Then we have
$$v(P_j)=\mbox{${1\over 2}$}(\ell^2+\ell+j),\qquad j=0, 1, \ldots , \ell-1,$$
$$v(H_{i,j})=\mbox{${1\over 2}$}(\ell^2+(2\ell-i)i+2(\ell-i)+j),\qquad j=0, 1, \ldots , \ell-2-i$$
and
$$v(H'_{i,j})=\mbox{${1\over 2}$}(\ell^2+(2\ell-i)i+3(\ell-i)-2+j),\qquad j=0, 1, \ldots , \ell-2-i.$$
It follows that the sequence $v(P_0)$, $v(P_1)$, $\ldots$ , $v(P_{\ell-1})$, $v(H_{0,0})$,
$v(H_{0,1})$, $\ldots$, $v(H_{0,\ell-2})$, $v(H'_{0,1})$, $\ldots $, $v(H'_{0,\ell-2})$,
$v(H_{1,0})$, $v(H_{1,1})$, $\ldots$, $v(H_{1,\ell-3})$, $v(H'_{1,1})$, $\ldots $,
$v(H'_{1,\ell-3})$, $v(H_{2,0})$, $\ldots$, $v(S_\ell)$ is exactly
the sequence $$\mbox{$1\over 2$}(\ell^2+k),\qquad k=\ell, \ell+1, \ldots , \ell^2.$$ This proves the assertion.
\hfill{$\square $}

\medskip\noindent{\bf Proof of Theorem 1.} The upper bound follows from (5) and the fact that
$$v^\ast (2,m)\le v(2,m).$$ Next, we prove the lower bound by modifying Arnold's ingenious method.

Let $\tau$ be a positive number and let $Q_\tau$ denote the set of all primitive integer vectors in
the domain $\{ (x,y):\ x^2+y^2\le \tau^2,\ x\le 0,\ y\ge 0\}$. One has
$$\sum_{{\bf v}\in Q_\tau }{\bf v}=(-\ell_\tau , \ell_\tau )$$
for a suitable integer $\ell_\tau$. Then,
we take ${\bf v}_0=(\ell_\tau , 0)$, ${\bf v}'_0=(0, -\ell_\tau )$ and define
$$V_\tau =\{ {\bf v}_0\} \cup Q_\tau \cup \{ {\bf v}'_0\}.$$
Let $M_\tau$ denote the convex lattice polygon whose oriented sides are all vectors in $V_\tau$
starting from ${\bf v}_0$ and ending with ${\bf v}'_0$ (see Figure 1).

$$\psannotate{\psboxto(13cm;0cm){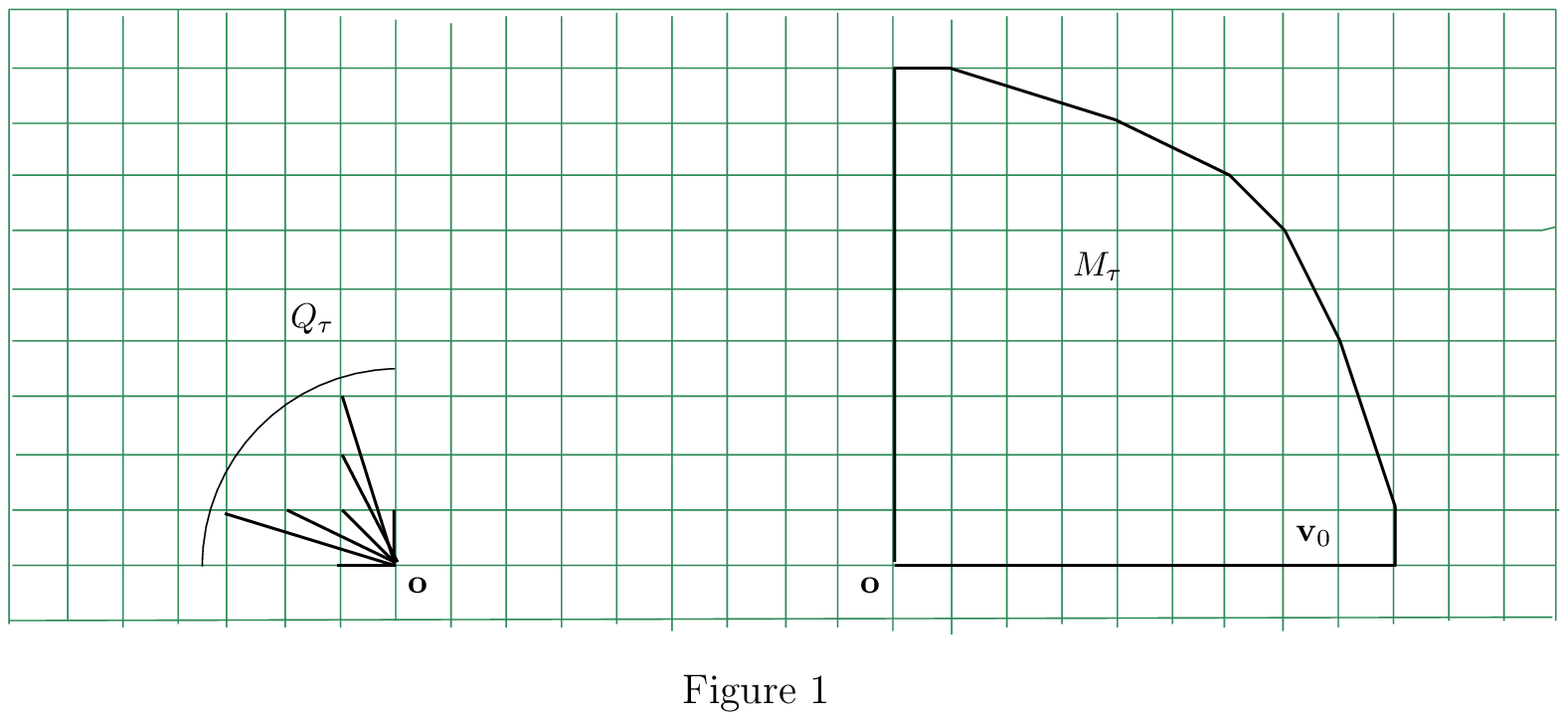}}{}$$

It is well-known in number theory (see \cite{hua82} p.125) that
$$|Q_\tau |={3\over {2\pi}}\tau^2+O(\tau\log \tau ).\eqno(8)$$
The convex lattice polygon $M_\tau$ has following properties:

\smallskip\noindent
{\bf 1.} It has $|V_\tau |$ vertices.

\smallskip\noindent
{\bf 2.} Let $C$ be the unit circular disc and let $$D=C\cap \{ (x,y):\ x, y\ge 0\}$$ be the
non-negative quadrant of $C$. Let $r$ be the largest number such that $rD\subseteq M_\tau$
and let $r'$ be the smallest number such that $M_\tau \subseteq r'D$. It can be easily deduced
from (8) that
$$\tau^3\ll r\le \ell_\tau \le r'\ll \tau^3.\eqno(9)$$

\smallskip\noindent
{\bf 3.} Each side, except the two longest, contains no other integral point except the ends.

\smallskip\noindent
{\bf 4.} The polygon $M_\tau$ changes only when $\tau^2$ passes through integral values, and
$${{v(M_{\tau +1})}\over {v(M_\tau )}}\le c$$
holds for some constant $c$. This inequality can be deduced from property 2.

\medskip
We will now construct polygons $P_i$ that are assembled with the help of smaller polygons
$P_i^1$, $P_i^2$, $P_i^3$ and $P_i^4$. For convenience, we enumerate the short sides of
$2M_\tau $ by $S_1$, $S_2$, $\ldots $,
$S_{|Q_\tau |}$ in anti-clockwise order. It follows by property 3 that each of these
sides (say $S_i$) has three integral points ${\bf p}_i^0$, ${\bf p}_i^1$ and ${\bf p}_i^2$
in this order. Then, ${\bf p}_1^0=(2\ell_\tau,0)$, ${\bf p}_i^0={\bf p}_{i-1}^2$ and
${\bf p}_{|Q_\tau |}^2=(0, 2\ell_\tau).$ For each vector $${\bf u}=(i_1, i_2, \ldots, i_{|Q_\tau |-1})\in
\{1, 2\}^{|Q_\tau |-1}$$ we obtain a convex lattice polygon
$$P_{\bf u}^1={\rm conv}\left\{ {\bf o}, {\bf p}_1^0, {\bf p}_1^{i_1},{\bf p}_2^{i_2}, \ldots ,
{\bf p}_{|Q_\tau |-1}^{i_{|Q_\tau |-1}}, {\bf p}_{|Q_\tau |}^2 \right\}$$
and, clearly, the lattice polygons $P_{\bf u}^1$, $P_{\bf v}^1$ differ for distinct ${\bf u}, {\bf v}\in \{ 1, 2\}^{|Q_\tau |-1}.$
So we obtain $2^{|Q_\tau |-1}$ convex lattice polygons. For sake of simplicity we enumerate them by
$P_1^1$, $P_2^1$, $\ldots $, $P_{2^{|Q_\tau |-1}}^1$ and
denote the set of these polygons by $\mathcal{P}^1$. It follows by property 3 that
$${1\over 2}\le {{v(P_i^1)}\over {v(2M_\tau )}}\le 1,\quad P_i^1\in \mathcal{P}^1.$$

Let $\bigtriangleup_\tau$ denote the lattice triangle with vertices $(0,0)$, $(0, 2\ell_\tau )$
and $(-2\ell_\tau , 0)$. Let $m$ be an even integer and choose $\tau $ to be the largest number
satisfying
$$v(2M_\tau )+v(\bigtriangleup_\tau )+5\ell_\tau \le {m\over 4}.\eqno (10)$$
Then
$$v(P_i^1)+2\ell_\tau^2+5\ell_\tau \le {m\over 4}$$
holds for all $P_i^1\in \mathcal{P}^1.$ Therefore, for each $P_i^1$ there are a positive integer $j$ and a corresponding number $\mu_i$ with $\ell_\tau \le \mu_i<5\ell_\tau $ satisfying
$${m\over 4}-v(P_i^1)-2\ell_\tau^2=4j\ell_\tau +\mu_i.\eqno (11)$$
By Lemma 2, one can extend
$\bigtriangleup_\tau $ to a lattice polygon $P^2_i$ contained in the square with vertices $(0,0)$,
$(0, 2\ell_\tau )$, $(-2\ell_\tau , 0)$ and $(-2\ell_\tau , 2\ell_\tau )$ and satisfying
$$v(P^2_i)=2\ell_\tau^2+\mu_i.\eqno (12)$$

Then we define
$$P^3_i=\{ (x,y):\ |x|\le 2\ell_\tau,\ 0\le y\le j\},$$
$$P^4_i=P^3_i\cup (P^2_i+(0,j))\cup (P^1_i+(0,j)),$$
$$P_i={\rm conv}\{ P^4_i\cup (- P^4_i)\}$$
and
$$\mathcal{P}=\{ P_1, P_2, \ldots , P_{2^{|Q_\tau |-1}}\}.$$

$$\hspace{-5cm}\psannotate{\psboxto(0cm;8cm){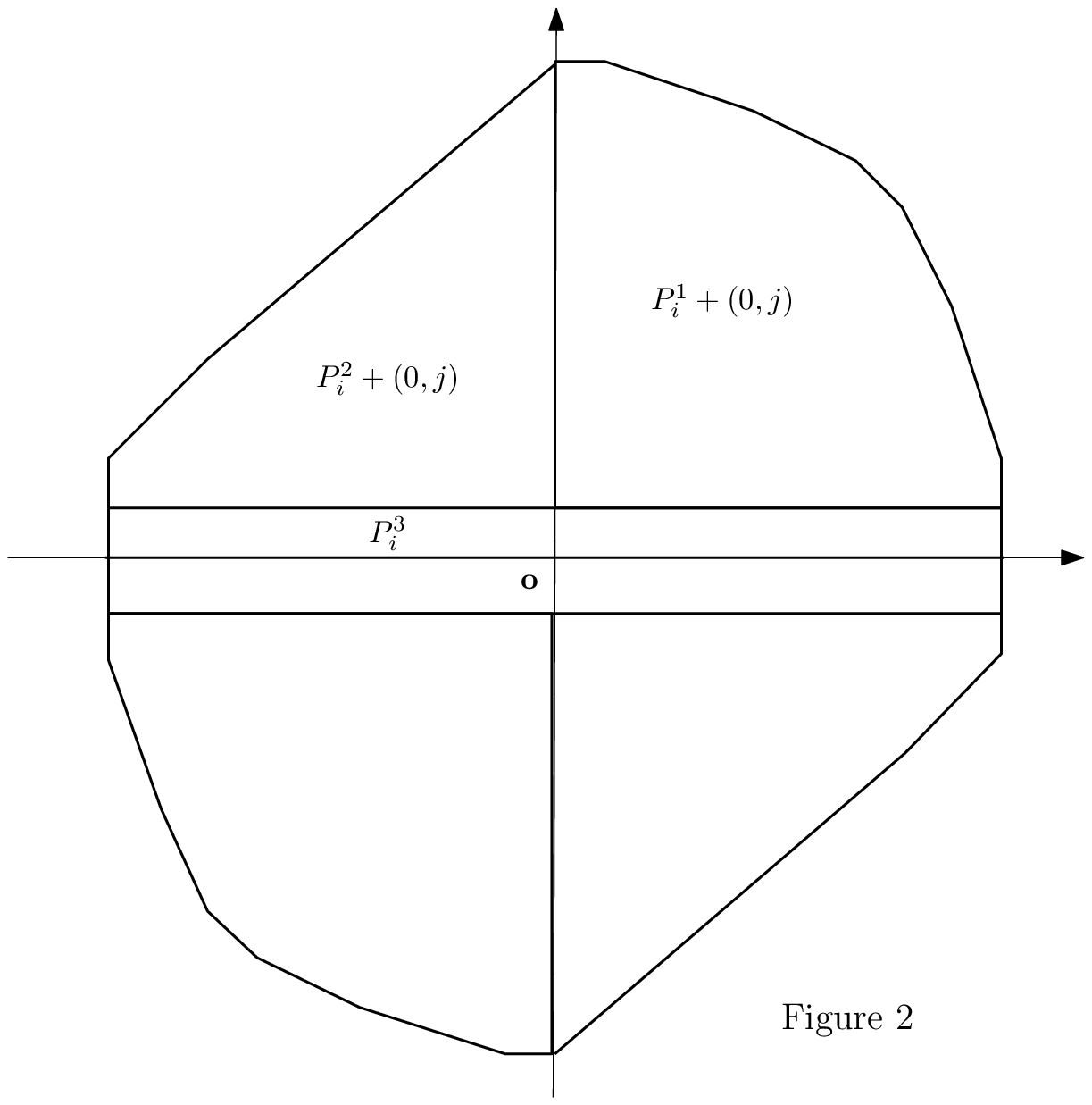}}{}$$

\noindent
Clearly, all $P_i$ are centrally symmetric convex lattice polygons as shown by Figure 2 and, by (11), (12) and their constructions,
$$v(P_i)={m\over 2},\qquad P_i\in \mathcal{P}.$$

Let $\ell (P_i)$ denote the maximal cardinality of the sets of collinear lattice points in $P_i$. It follows
by the constructions of $P_i$ that
$$\ell (P_i)=2(j+2\ell_\tau )+1$$
and the lattice segment from $(0,j+2\ell_\tau )$ to $(0, -j-2\ell_\tau )$ is the only longest one passing
the origin. Thus, any pair of the polygons in $\mathcal{P}$ are not equivalent and hence
$$v^\ast (2,m)\ge |\mathcal{P}|=2^{|Q_\tau |-1}.\eqno (13)$$

On the other hand, by the maximum assumption on $\tau $ in (10) and (9), we have
$$m\le 4(v(2M_{\tau +1})+2\ell_{\tau +1}^2+5\ell_{\tau +1} )\ll \ell_{\tau +1}^2 \ll \tau^6.\eqno(14)$$
Thus, by (13), (8) and (14) we get
$$\log v^\ast (2,m)\gg |Q_\tau |\gg \tau^2\gg m^{1\over 3}.$$
This concludes the proof of the theorem. \hfill{$\square $}

\medskip
\noindent
{\bf Remark 2.} B\'ar\'any \cite{bara08'} proposed the following problem: Decide whether
$$\lim_{m\to\infty}m^{-{1\over 3}}\log v(2, m)$$
exists or not, and determine the limit if it exists. Of course, one can ask the same question for $v^*(2,m)$.

\medskip
Let $F(d,m)$ denote the family of all $d$-dimensional convex lattice
polytopes of volume $m/d!$. The family is divided by the relation
$\sim $ into $v(d,m)$ classes. By choosing one from each class we get a
representative set for $F(d,m)$. To determine a representative set
for $F(d,m)$ is an interesting problem as well. For this purpose,
we will introduce some invariants in next sections.

In the plane, when $m$ is small, based on Lemma 1 and Pick's theorem
we can determine representative sets for $F(2,m)$ (see Appendix 1) and therefore
the values of $v(2,m)$. The methodology will be introduced in Section 4.

\medskip
\centerline{
\renewcommand\arraystretch{1.5}
\begin{tabular}{|c|c|c|c|c|c|c|c|}
\hline $m$ & $1$ & $2$ &
$3$ & $4$ & $5$ & $6$ & $7$\\
\hline $v(2,m)$ & $1$ & $2$ & $3$ & $7$ & $6$ & $13$ & $11$\\
\hline
\end{tabular}}

\section*{4. Convex Lattice Polygons of Given Cardinality}

In this section we study the classification problem for convex lattice polygons of
given cardinality. First, we prove the following result which is an analogue to the results of
Arnold, B\'ar\'any and Pach.

\medskip\noindent
{\bf Theorem 2.} {\it When $w$ is sufficiently large, we have}
$$w^{1\over 3}\ll \log \kappa (2, w)\ll w^{1\over 3}.$$

\medskip\noindent{\bf Proof.} If $|P|=w$ and $v(P)=m/2$ hold for suitable positive integers
$w$ and $m$, by Pick's theorem we get
$$w-2\le m\le 2w-2.$$ Then by (5) we have
\begin{eqnarray*}
\log \kappa (2,w) &\le & \log \left( \sum_{m=w-2}^{2w-2}v(2,m)\right)\\
&\le & \log \left((w+1)\cdot \max_{w-2\le m\le 2w-2}v(2, 2w)\right)\\
&\ll & w^{1\over 3}+\log w\\
&\ll & w^{1\over 3},
\end{eqnarray*}
which proves the upper bound.

Next, we prove the lower bound by following Arnold's process.

Let $\tau$ be a large number, let $V_\tau$ denote the set of all primitive integer vectors in the
semicircle $\{ (x,y):\ x^2+y^2\le \tau^2,\ x>0\}$, and let $M_\tau$ denote the convex lattice
polygon whose oriented sides are all vectors in $V_\tau$ and $-\sum_{{\bf v}\in V_\tau}{\bf v}$ (see Figure 3).

$$\psannotate{\psboxto(13cm;0cm){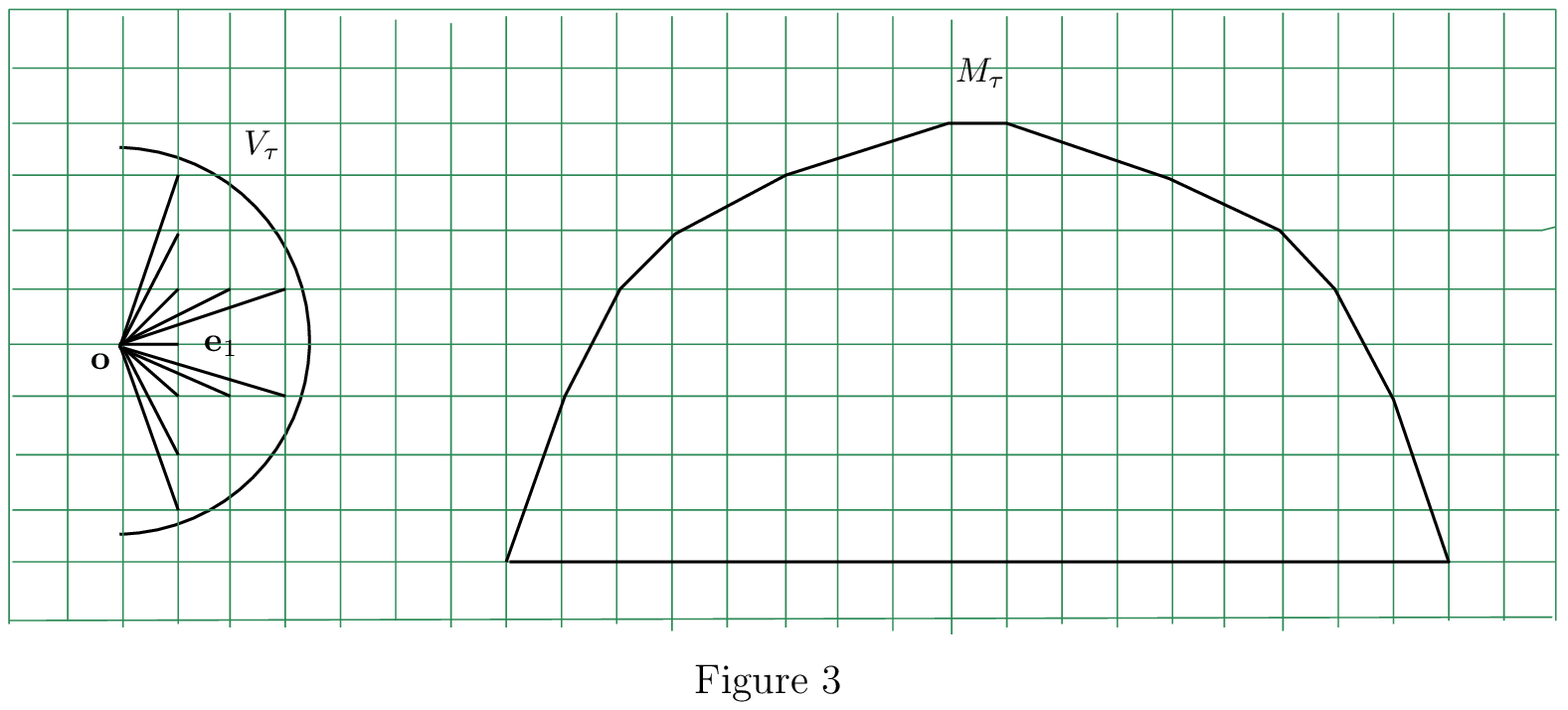}}{}$$

From (8) we directly deduce
$$|V_\tau |={3\over \pi}\tau^2+O(\tau\log \tau ).\eqno(15)$$
The convex lattice polygon $M_\tau$ has following properties:

\noindent
{\bf 1.} It has $|V_\tau |+1$ vertices.

\noindent
{\bf 2.} Let $r$ be the largest radius of semicircular discs contained in $M_\tau$ and let $r'$ be the
smallest radius of semicircular discs containing $M_\tau $. By (15) it can be easily deduced that
$$\tau^3\ll r\le r'\ll \tau^3.\eqno(16)$$

\noindent
{\bf 3.} Each side, except the diameter, contains no other integral point except the ends. The
oriented edge at the top is ${\bf e}_1$ (see Figure 3).

\noindent
{\bf 4.} The polygon $M_\tau$ changes only when $\tau^2$ passes through integral values, and
$${{|M_{\tau +1}|}\over {|M_\tau |}}\le c_1\eqno(17)$$
holds for some constant $c_1$. This inequality can be deduced from property 2.

We will now construct polygons $P_i$ that are assembled with the help of smaller polygons
$P_i^1$ and $P_i^2$. For convenience, we enumerate the short sides of $2M_\tau $ by $S_1$, $S_2$, $\ldots $,
$S_{|V_\tau |}$ in the clockwise order. It follows by property 3 that each of these sides (say $S_i$)
has three integral points ${\bf p}_i^0$, ${\bf p}_i^1$ and ${\bf p}_i^2$ in the order. Then,
${\bf p}_1^0=(0,0)$, ${\bf p}_i^0={\bf p}_{i-1}^2$ and ${\bf p}_{|V_\tau |}^2=2\sum_{{\bf v}
\in V_\tau }{\bf v}.$ For each vector $${\bf u}=(i_1, i_2, \ldots, i_{|V_\tau |-1})\in
\{1, 2\}^{|V_\tau |-1}$$ we obtain a convex lattice polygon
$$P_{\bf u}^1={\rm conv}\left\{ {\bf p}_1^0, {\bf p}_1^{i_1},{\bf p}_2^{i_2}, \ldots , {\bf p}_{|V_\tau |-1}^{i_{|V_\tau |-1}}, {\bf p}_{|V_\tau |}^2\right\}$$
and, clearly, the lattice polygons $P_{\bf u}^1$, $P_{\bf v}^1$ differ for distinct ${\bf u},{\bf v}\in \{ 1, 2\}^{|V_\tau |-1}.$ So we obtain $2^{|V_\tau |-1}$ convex lattice polygons. We enumerate them by
$P_1^1$, $P_2^1$, $\ldots $, $P_{2^{|V_\tau |-1}}^1$ and denote
the set of these polygons by $\mathcal{P}^1$. It follows by property 3 that
$${1\over 2}\le {{|P_i^1|}\over {|2M_\tau |}}\le 1,\quad P_i^1\in \mathcal{P}^1.$$

Assume that $w$ is an integer and $\tau$ is the largest number satisfying
$$|2M_\tau |\le w.\eqno (18)$$
Then we have
$$|P_i^1|\le w, \quad P_i^1\in \mathcal{P}^1.$$
Let $n$ be an integer determined by
$$2\sum_{{\bf v}\in V_\tau }{\bf v}=(n,0).$$
There are nonnegative integers $k$ and $\ell$ with $\ell\le n$ satisfying
$$w-|P_i^1|=k(n+1)+\ell.$$
Then we define
$$P^2_i=\{ (x,y):\ 0\le x\le n,\ -k\le y\le -1;\ y=-k-1,\ 0\le x\le \ell -1\},$$
$$P_i={\rm conv}\{ P_i^1\cup P^2_i\},$$
and define $\mathcal{P}$ to be the set of all these polygons $P_i$. Clearly, all $P_i$ are convex
lattice polygons with
$$|P_i|=w,$$
and for each $P_i$ there is at most another $P_j\in \mathcal{P}$ satisfying $P_i\sim P_j$.
Therefore, we get
$$\kappa (2, w)\ge {{|\mathcal{P}|}\over 2}=2^{|V_\tau |-2}.\eqno (19)$$
On the other hand, by the maximum assumption on $\tau $ in (18), (17) and (16) we have
$$w\le |2M_{\tau +1}|\ll |M_\tau |\ll \tau^3\cdot \tau^3=\tau^6.\eqno(20)$$
Thus, by (19), (15) and (20) we get
$$\log \kappa (2,w)\gg |V_\tau |\gg \tau^2\gg w^{1\over 3}.$$
The theorem is proved. \hfill{$\square $}

\medskip\noindent
{\bf Theorem 3.} {\it When $w$ is odd and sufficiently large, we have}
$$w^{1\over 3}\ll \log \kappa^* (2, w)\ll w^{1\over 3}.$$

\noindent
{\bf Remark 3.} The upper bound of Theorem 3 follows from Theorem 2 and the fact that
$$\kappa^*(2,w)\le \kappa (2, w).$$
The lower bound can be proved by modifying the proof of Theorem 1, simply replacing
the areas by cardinalities. In fact, Lemma 2 is no longer needed now.

\medskip
Similar to B\'ar\'any's problem, the following one seems interesting and challenging as well.

\noindent
{\bf Problem 1.} Decide whether
$$\lim_{w\to\infty}w^{-{1\over 3}}\log \kappa (2, w)$$ and
$$\lim_{w\to\infty}w^{-{1\over 3}}\log \kappa^\ast (2, w)$$
exist or not. Determine the limits if they exist.

\medskip
Let $G(d,w)$ denote the family of all $d$-dimensional convex lattice
polytopes of cardinality $w$. Clearly the family is divided by the relation
$\sim $ into $\kappa (d,w)$ classes. By choosing one from each class we get a
representative set for $G(d,w)$.

\medskip
To obtain a representative set for $G(2,w)$ (as well as for $F(2,m)$, which was defined
directly after Remark 2) for a given positive integer $w$, we follow the following steps:

\medskip\noindent
{\bf Step 1.} Let $\ell (P)$ denote the maximum number such that $P$ has $\ell (P)$ collinear integer
points.  Let $\lfloor x \rfloor$ denote the largest integer $z$ satisfying $z\le x$ and
let $\lceil x\rceil$ denote the smallest integer $z$ satisfying $x\le z$. According to Lemma 1, we have
$$\ell (P)\ge \lceil \sqrt{w} \rceil , \qquad P\in G(2, w).$$
Let $L(P)$ be such a maximal collinear set. It is well-known in geometry of numbers that there is
a $\mathbb{Z}^2$-preserving affine transformation which transfers $L(P)$ to $\ell (P)$ successive integer points on the $x$-axis. Therefore, without loss of generality, we assume that
$$L(P)=\left\{ (-\lfloor \mbox{$1\over 2$}(\ell (P)-1)\rfloor , 0), (-\lfloor \mbox{$1\over 2$}(\ell (P)-1)
\rfloor +1, 0), \ldots, (\lceil \mbox{$1\over 2$}(\ell (P)-1)\rceil , 0)\right\} .$$

\medskip\noindent
{\bf Step 2.} We further classify the convex lattice polygons according to the number of points above and
below the $x$-axis. Let $\ell (P)$ take the values $w-1$, $w-2$, $\ldots$, $\lceil \sqrt{w}\rceil $, respectively.
For each value, we get a list of partitions $$w-\ell (P)=z_1+z_2,$$ where $z_1$ and $z_2$ are nonnegative
integers with $z_1\ge z_2$.

Let $L_i(P)$ denote the set of $P\cap \{ (x,y):\ y=i, x\in \mathbb{Z}\}$ and let $\ell_i(P)$ denote the cardinality of $L_i(P)$. By a suitable $\mathbb{Z}^2$-preserving affine transformation we can take
$$L_1(P)=\left\{ (-\lfloor \mbox{$1\over 2$}(\ell_1(P)-1)\rfloor , 1), (-\lfloor \mbox{$1\over 2$}(\ell_1(P)-1)
\rfloor +1, 1), \ldots, (\lceil \mbox{$1\over 2$}(\ell_1(P)-1)\rceil , 1)\right\} .$$ In particular, $(0,1)\in L_1(P)$. For convenience, we abbreviate $P\cap \{ (x,y):\ y\ge 0\}$ and $P\cap \{ (x,y):\ y\le 0\}$ by $P'$ and $P^*$, respectively.

Let $h$ be the maximal integer that $\ell_h(P)\not= 0$, $F$ be the area of $P'$, $I$ be the number of the interior lattice points of $P'$ and let $J$ be the number of the lattice points on the boundary of $P'$. By Pick's theorem, we get
$$h\cdot \ell (P)\le  2F=J+2I-2\le 2(J+I)=2(\ell (P)+z_1)$$
and therefore
$$h\le 2+{{z_1}\over {\ell (P)}}.$$

By convexity, $P'\cap \mathbb{Z}^2$ is contained in the region bounded by four lines $\{ (x,y):\ y=0\}$, $\{(x,y):\ y=x-\lceil {1\over 2} (\ell (P)-1)\rceil \}$, $\{(x,y):\ y=2+\lfloor {{z_1}\over {\ell (P)}}\rfloor \}$ and $\{(x,y):\ y=-x-\lfloor {1\over 2} (\ell (P)-1)\rfloor \}$. Similarly, the lower part $P^*\cap \mathbb{Z}^2$ is contained in the region bounded by $\{(x,y):\ y=0\}$, the line passing $(0,1)$ and $(\lceil {1\over 2}(\ell (P)-1)\rceil +1, 0)$, $\{ (x,y):\ y=-2-\lfloor {{z_2}\over {\ell (P)}}\rfloor\}$, and the line passing $(0,1)$ and $(-\lfloor {1\over 2}(\ell (P)-1)\rfloor -1, 0)$.

Thus, for each partition we can routinely get all possible convex lattice polygons (in the sense of the equivalence) with $z_1$ lattice points above $L(P)$ and $z_2$ points below $L(P)$. Thus we get a set $\mathcal{F}$ which contains a representative set of $G(2,w)$.

\medskip\noindent
{\bf Step 3.} To determine two polygons are equivalent or not, we need to define a set of invariants under $\sim$.
As an example, for a convex lattice polygon $P_i\in \mathcal{F}$, we define an invariant vector ${\bf v}_i=(v_i^1, v_i^2, v_i^3),$ where $v_i^1=\ell (P_i)$, $v_i^2$ is the number of vertices of $P_i$, and $v_i^3$ is the number of the interior lattice points of $P_i$. Clearly, we have ${\bf v}_i={\bf v}_j$ if $P_i\sim P_j$. Therefore, $P_i\not \sim P_j$ if ${\bf v}_i\not={\bf v}_j$.

When $w$ is small, a couple of nice invariants are enough to distinguish all non-equivalent lattice polygons. For large $w$, it seems that a complicated set of invariants is required.

Following these steps and applying the invariant vector introduced above, we get representative sets for $G(2,3)$, $G(2,4)$, $G(2,5)$, $G(2,6)$ and $G(2,7)$,
as listed in Appendix 2. Consequently, we get the exact values of $\kappa (2,w)$ for $3\le w\le 7$. Up to now we have not employed a computer in this project. It is possible to create a computer program based on these steps to determine $G(2,w)$ and $\kappa (2, w)$ for some large $w$.

\medskip
\centerline{
\renewcommand\arraystretch{1.5}
\begin{tabular}{|c|c|c|c|c|c|}
\hline $w$ & $3$ & $4$ & $5$ & $6$ & $7$\\
\hline $\kappa (2,w)$ & $1$ & $3$ & $6$ & $13$ & $21$ \\
\hline
\end{tabular}}

\section*{5. Convex Lattice Polytopes of Given Cardinality}

In this section we study $\kappa (d, w)$ and $\kappa^*(d,w)$ for $d\ge 3$.

Let $\{ {\bf e}_1, {\bf e}_2, \ldots , {\bf e}_d\}$ be an orthonormal base of $\mathbb{E}^d$. We define $P(d,w,k)$ to be the convex hull of
$$\left\{ {\bf o}, {\bf e}_i,\ -j{\bf e}_d,\  \sum_{n=1}^{d-1}{\bf e}_n +k{\bf e}_d:\ i=1, 2,
\ldots, d-1;\ j=1, \ldots , w-d-1 \right\}.$$
When $d\ge 3$ and $w\ge d+1$, it can be easily shown that $P(d,w,k)$ is a $d$-dimensional convex lattice polytope with
$$|P(d,w,k)|=w.\eqno (21)$$
For example, Figure 4 shows the three-dimensional convex lattice polytope $P(3,5,3)$.
$$\psannotate{\psboxto(9.6cm;0cm){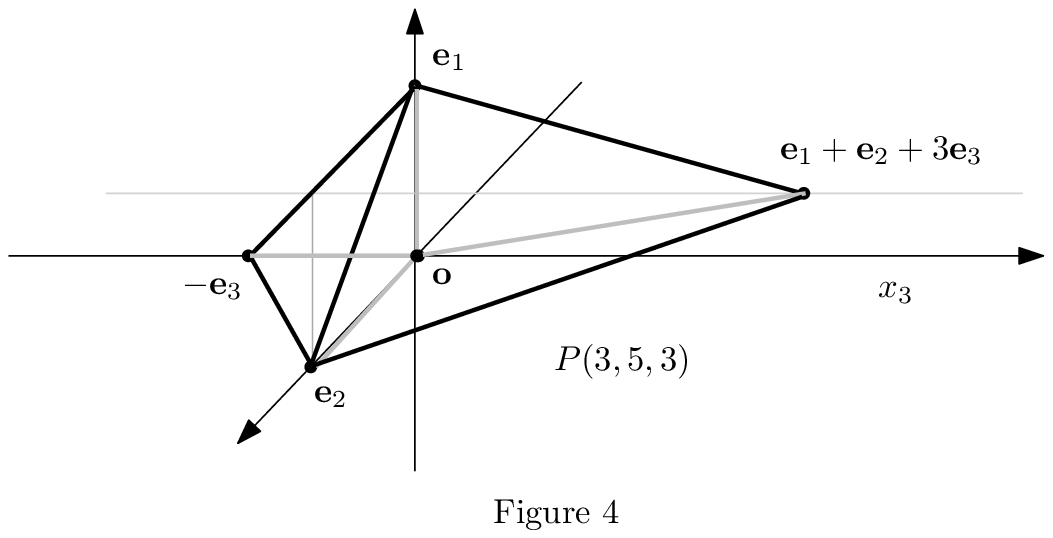}}{}$$

When $k\ge (d-2)(w-d-1)$, the lattice polytope $P(d,w,k)$ can be divided into two simplices
$$P_1={\rm conv} \left\{ {\bf o}, {\bf e}_i,\  \sum_{n=1}^{d-1}{\bf e}_n +k{\bf e}_d:\ i=1, 2,
\ldots, d-1 \right\}$$
and
$$P_2={\rm conv}\left\{ {\bf o}, {\bf e}_i,\ -j{\bf e}_d,:\ i=1, 2,
\ldots, d-1;\ j=1, \ldots , w-d-1 \right\}.$$
Thus, we get
$$v(P(d,w,k))={{k+w-d-1}\over {d!}}.\eqno (22)$$

It follows by (21) that all $P(d,w,k)$ belong to $G(d,w)$. On the other hand, it follows by (22) that
$$P(d,w,k_1)\not\sim P(d,w,k_2)$$
whenever $k_i\ge (d-2)(w-d-1)$ for $i=1$ and $2$. Thus we have proved the following result.

\medskip\noindent
{\bf Theorem 4.} {\it When $d\ge 3$ and $w\ge d+1$, we have}
$$\kappa (d, w)=\infty .$$

\medskip\noindent
{\bf Remark 4.} Comparing with (2), one can see the essential difference between $v(d,m)$ and
$\kappa (d, w)$ when $d\ge 3$: the first is finite, the second is infinite.

\medskip
Let $i(P)$ denote the number of the interior lattice points in $P$. Improving a result of Lagarias and
Ziegler \cite{laga91}, O. Pikhurko in 2001 obtained a upper bound (see (9) of \cite{pikh01}) which implies that
$$v(P)\le c_d\cdot |P|\eqno (23)$$
if $i(P)\not= 0$, where $c_d$ is a constant depends only on $d$. Let $\kappa '(d,w)$ denote the number of different classes of $d$-dimensional convex lattice polytopes $P$ with $|P|=w$ and $i(P)\not= 0$. It follows by (23) and (2) that

\begin{eqnarray*}
\log \kappa ' (d,w)&\le &\log \left(\sum_{m=1}^{c_d\cdot w}v(d,m)\right)\\
&\le &\log \left( c_d\cdot w \max_{m\in \{ 1, 2, \ldots, c_d\cdot w\}}\ \{v(d, m)\}\right)\\
&= &\log \left( c_d\cdot w\right)+ \max_{m\in \{ 1, 2, \ldots, c_d\cdot w\}}\ \{\log v(d, m)\}\\
&\ll & w^{{d-1}\over {d+1}}.
\end{eqnarray*}
In particular, we have $i(P)\not= 0$ if $P$ is centrally symmetric. Thus, we get
$$\log \kappa^*(d,w)\le \log \kappa '(d,w).$$ 
As a conclusion, we obtain the following result.

\bigskip\noindent
{\bf Theorem 5.} {\it We have
$$\log \kappa '(d, w)\ll w^{{d-1}\over {d+1}}$$
and}
$$\log \kappa^\ast (d, w)\ll w^{{d-1}\over {d+1}}.$$

\medskip
\noindent
{\bf Remark 5.} Comparing with Theorem 4, it shows the essential difference between $\kappa (d,w)$ and $\kappa^\ast (d,w)$ when $d\ge 3$.

\bigskip\bigskip\noindent
{\large\bf Acknowledgement.} We are grateful to Professor Imre B\'ar\'any, Matthias Henze and the referee for their helpful comments and remarks. The quality of this paper has been much improved by the suggestions of the referee and Matthias Henze.

\bigskip
\bibliographystyle{amsplain}

\bigskip\bigskip
\noindent H. Liu and C. Zong, School of Mathematical Sciences, Peking
University, Beijing 100871, China

\noindent E-mail: cmzong@math.pku.edu.cn

\vfill\eject
\section*{Appendix 1: Representative Sets for $F(2,m)$}

\medskip
$$\psannotate{\psboxto(12cm;0cm){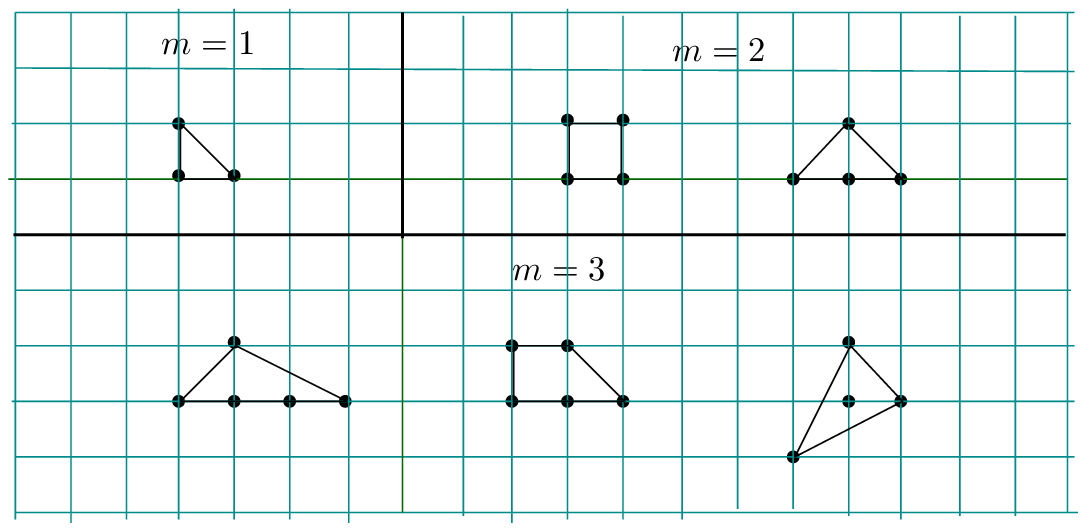}}{}$$

$$\psannotate{\psboxto(12cm;0cm){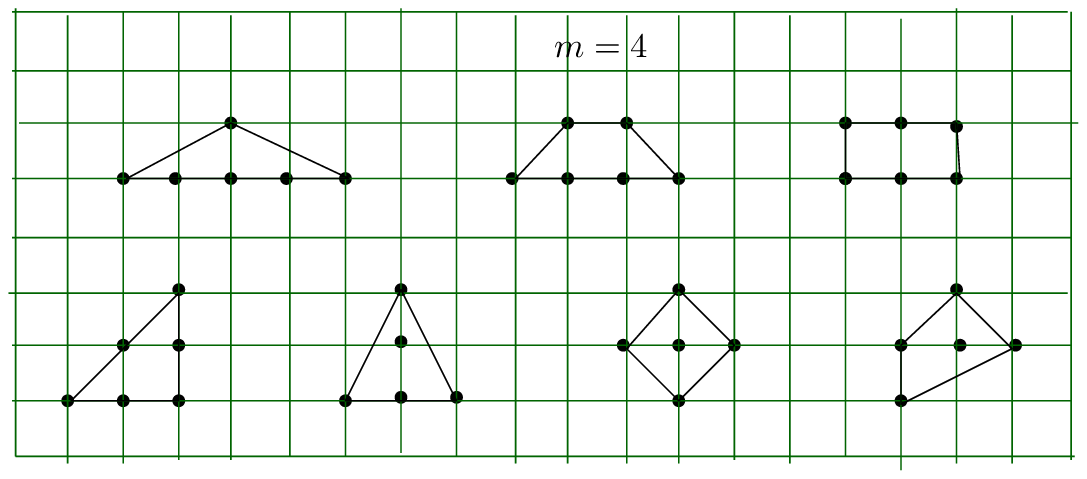}}{}$$

$$\psannotate{\psboxto(12cm;0cm){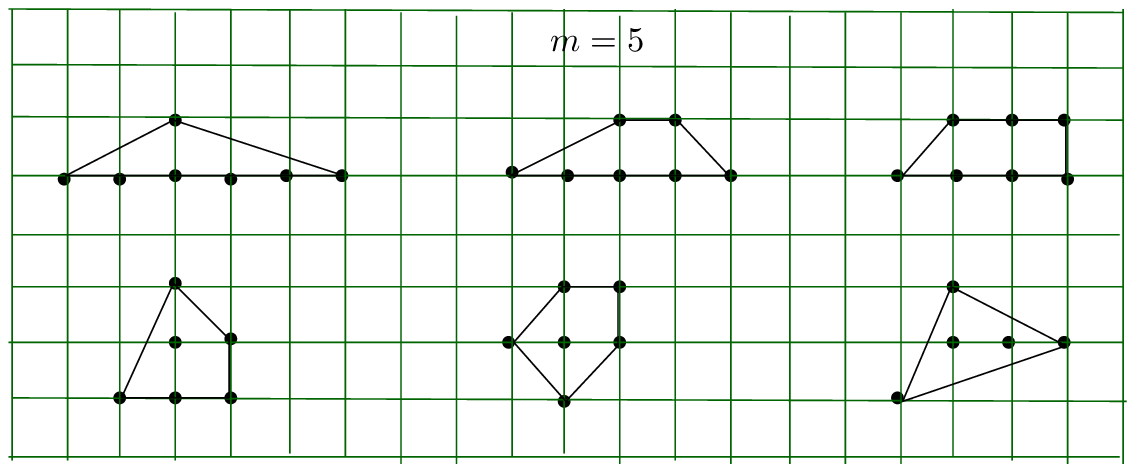}}{}$$

$$\psannotate{\psboxto(12cm;0cm){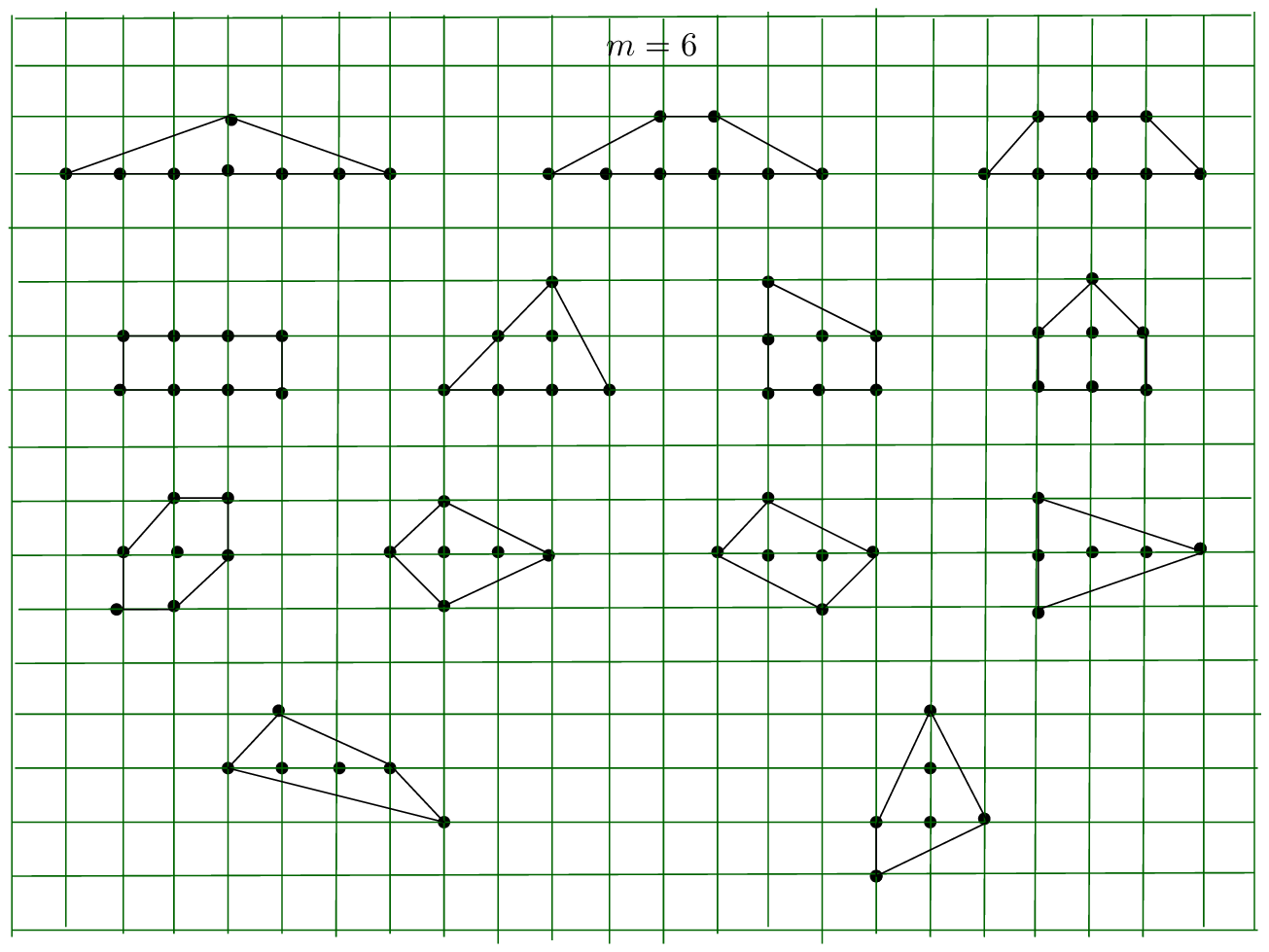}}{}$$

$$\psannotate{\psboxto(12cm;0cm){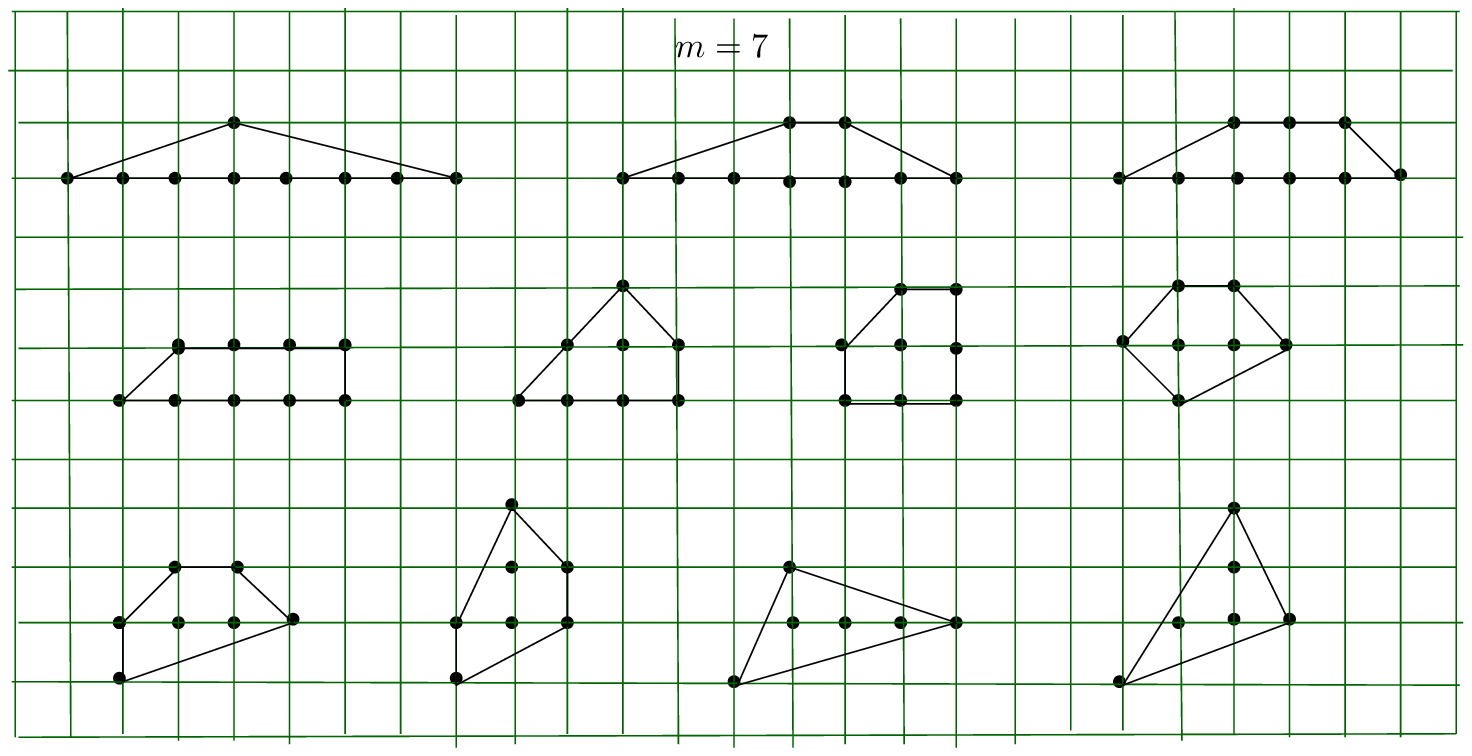}}{}$$

\section*{Appendix 2: Representative Sets for $G(2,w)$}

\medskip
$$\psannotate{\psboxto(12cm;0cm){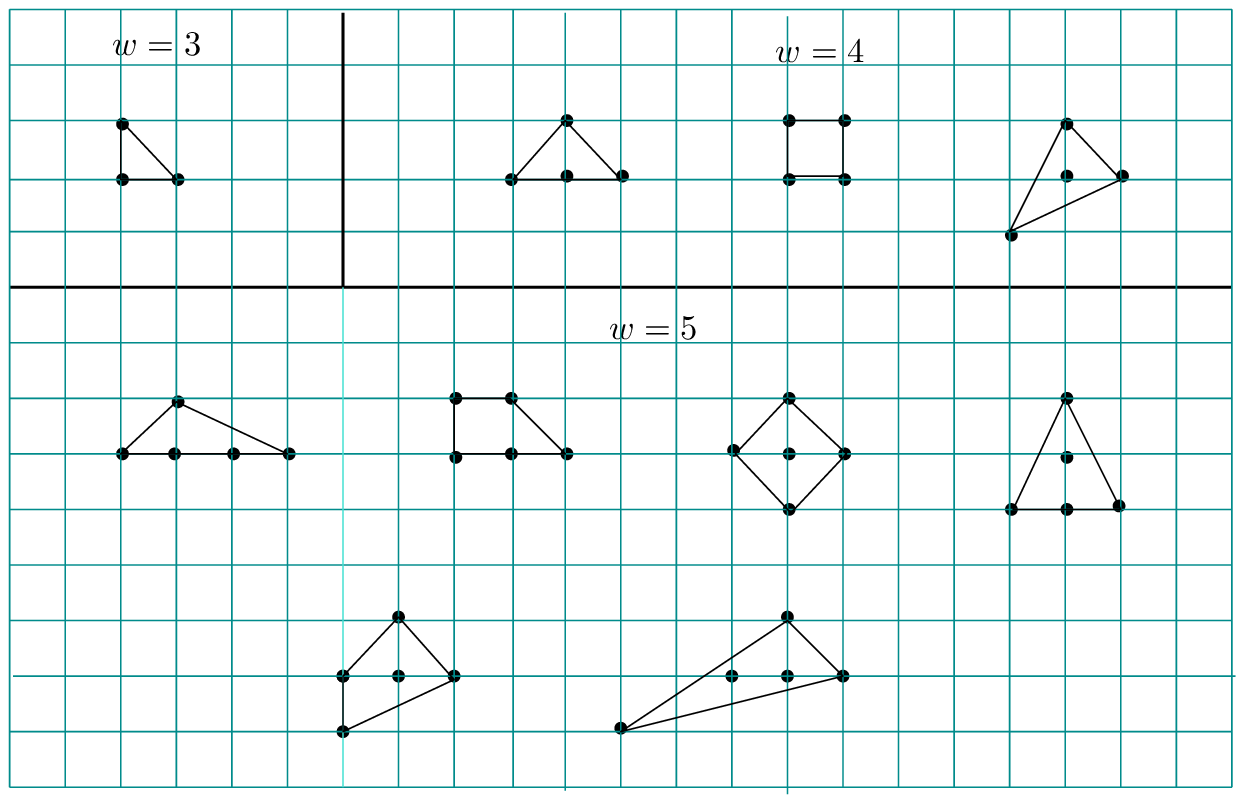}}{}$$

$$\psannotate{\psboxto(12cm;0cm){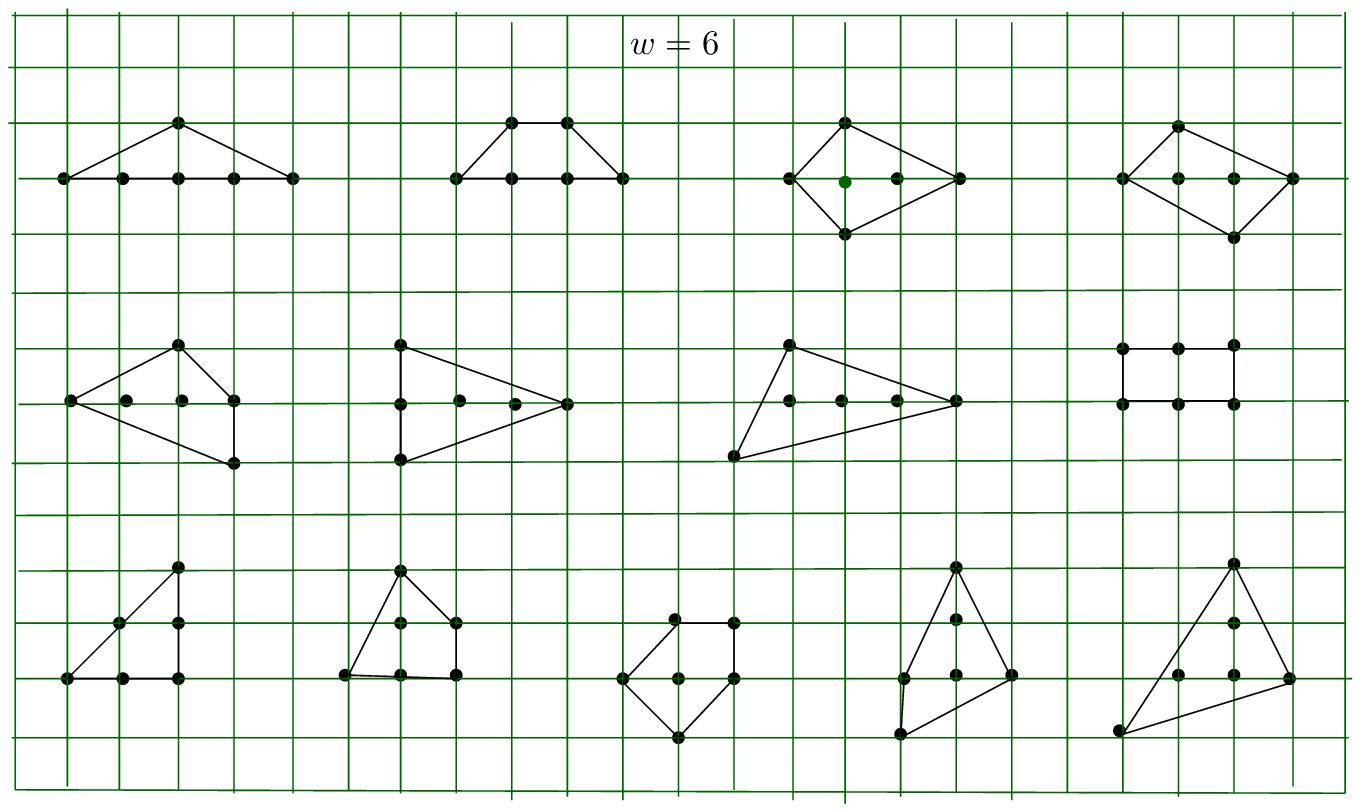}}{}$$

$$\psannotate{\psboxto(12cm;0cm){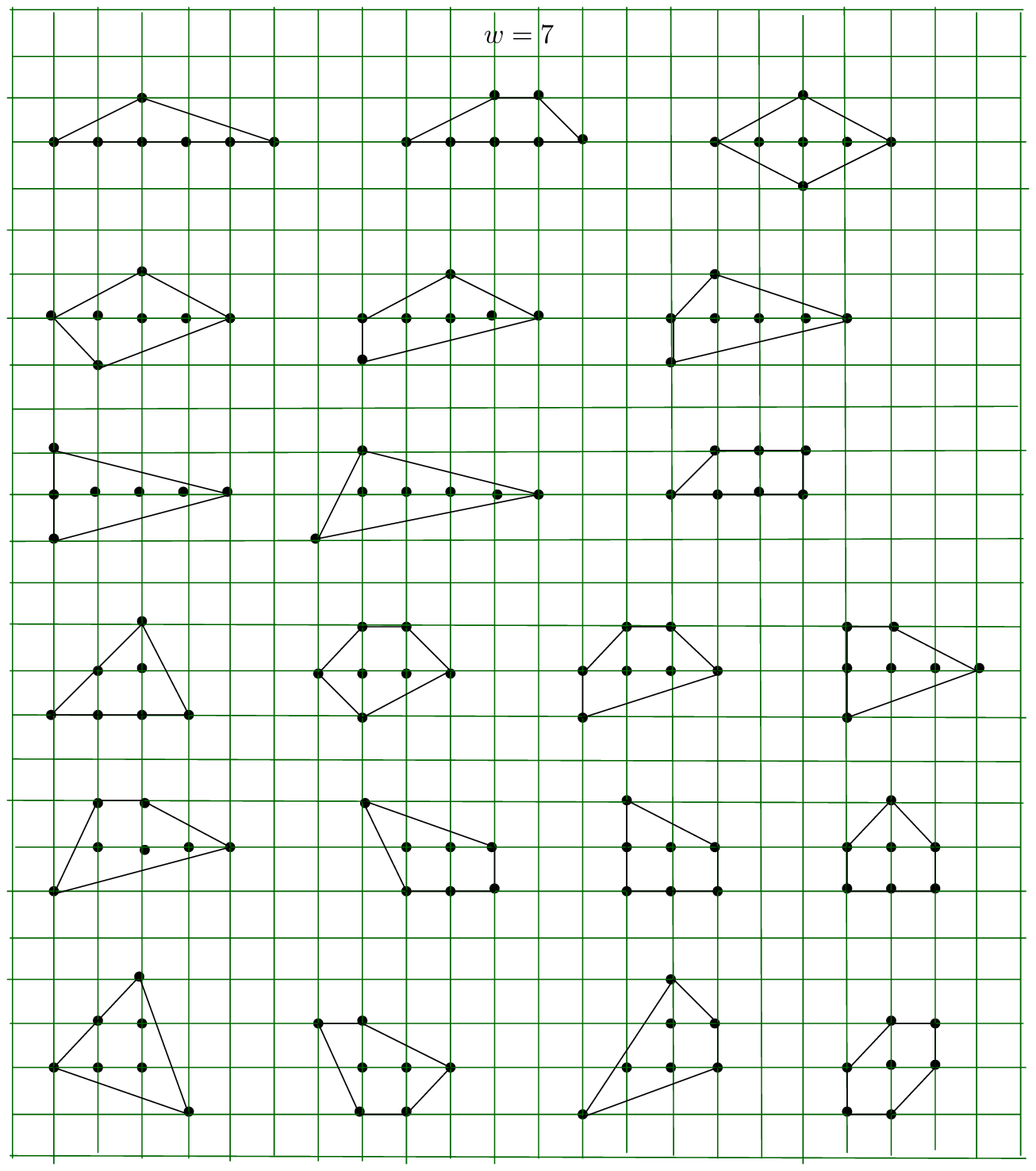}}{}$$

\end{document}